\documentclass[twoside,11pt]{article}   	
\usepackage{geometry}                		
\usepackage{graphicx}				
\usepackage{mathtools}
\usepackage{mathrsfs}
\usepackage[all]{xy}
\usepackage{amsmath}
\usepackage{MnSymbol}

\usepackage{wasysym}
\usepackage{tikz,tikz-cd}
\newtheorem{theorem}{Theorem}
\newtheorem{lemma}{Lemma}
\newtheorem{corollary}{Corollary}

\begin{document} 
 \title {Lefschetz theorem for holomorphic one-forms on weakly 1-complete manifolds}

\author{Chen Zhou}
\date{}
\maketitle
\begin{abstract}
  For a holomorphic one-form $\mathbf{\xi}$ on a weakly 1-complete manifold $X$ with certain properties, we will discuss the connectivity of the pair $(\hat{X},F^{-1}(z))$, where $\pi:\hat{X} \to X$ is a covering map and $F$ is a holomorphic function on $\hat{X}$ such that $dF=\pi^*\mathbf{\xi}$. We will also discuss the criteria about when such a manifold X admits a proper holomorphic mapping onto a Riemann surface.
	
\end{abstract}
\section{Introduction}	
 Simpson considered the Lefschetz theorem for a holomorphic one-form on a  compact complex manifold in \cite{r11}. He showed that if $X$ is an algebraic manifold with a holomorphic one-form $\mathbf{\xi}$ on it and if $\pi:\hat{X}\to X$ is a covering map such that the pullback of $\mathbf{\xi}$ is exact (that is, $\pi^*\mathbf{\xi}=dF$ for some holomorphic function $F$ on $\hat{X}$), then theorems about the connectivity of the pair $(\hat{X}, F^{-1}(z))$ can be obtained. He assumed $X$ is an algebraic manifold there,  but his method works for any compact complex manifold. 

The idea Simpson used to prove the theorems about the connectivity of the pair $(\hat{X}, F^{-1}(z))$ comes from Morse theory. As pointed in \cite{r11}, $\hat{X}$ is obtained from a fixed fiber $F^{-1}(z)$ with other fibers attached at the singular points of $F$. In order to study the connectivity of the pair $(\hat{X}, F^{-1}(z))$, we need to analyze the local topology in the neighborhood of every singular fiber, and we need the neighborhoods of two different singular fibers not too near to each other. This is the case when the manifold $X$ is compact, since then the subspace $A$ of $X$ where $\mathbf{\xi}$ vanishes has only finitely many components, and we can construct neighborhoods for these components (far away enough from each other) and then lift these neighborhoods to $\hat{X}$. However, when the manifold $X$ is not compact, such neighborhoods may not exist. But we will see below that in some special situations, we can still employ Simpson's idea, and after some modification, to get a similar result. In \cite{r11} Simpson also showed that if the pair $(\hat{X}, F^{-1}(z))$ is not 1-connected, then there exists a holomorphic mapping from $X$ to a Riemann surface, and $\mathbf{\xi}$ is the pullback of a one-form on that Riemann surface. We will see this remains true for some certain class of weakly 1-complete manifolds (Recall that a complex manifold is called a weakly 1-complete manifold if it admits a continuous plurisubharmonic exhaustive function).

   As pointed out in \cite{r11}, this connectivity question arises when we consider about equivariant harmonic maps from the universal covering of $X$ to trees (see \cite{r20}); the fibers of the harmonic map are unions of connected components of the integral leaves of a harmonic one-form, so it is essential to discuss the connectivity of these leaves. And we can also see that in some situations, the leaf space has the structure of a Riemann surface, making the quotient map from $X$ holomorphic (see \textbf{Proposition 3.5} in \cite{r22}). 
   
    In case of K\"ahler manifolds (which need not to be compact), it is known that harmonic maps from them to trees are pluriharmonic under some assumptions (see \cite{r18}), since trees are special kind of negatively curved spaces. An interesting  case now is when the pluriharmonic map $u$ from a non-compact K\"ahler manifold $X$ to a tree is proper. Then $u$ gives rise to a holomorphic quadratic differential on $X$ (See \cite{r18}, \cite{r20}). We know there exists a branched double covering $\Pi:\breve{X}\to X $ so that the pull-back of the quadratic differential by $\Pi$ equals the square of some  holomorphic one-form $\mathbf{\xi}$ on $\breve{X}$. Such $\mathbf{\xi}$ equals $\partial\varphi$ for some proper pluriharmonic function $\varphi$ on $(\breve{X}-K)$, where $K$ is a compact subset of $\breve{X}$.
   
   This motivated us to consider the Lefschetz type theorems on a class of K\"ahler manifolds which admit such $\varphi$, since then the properness of $\varphi$ enables us to construct the neighborhoods (far away enough from each other) of different components of $A$ (recall $A$ is the subspace of $X$ where $\mathbf{\xi}$ vanishes) lying inside some $\varphi^{-1}[a,b]$. So now we can modify the idea in \cite{r11} to get the following theorem,

\begin{theorem}\label{thm1}
	Let $X$ be a $d$-dimensional $(d\geqslant2)$ non-compact connected K\"ahler manifold, $\mathbf{\xi}$ a holomorphic one-form on $X$ which is not identically zero. Suppose that $\mathbf{\xi}=\partial\varphi$ on $(X-K)$ for some compact subset $K$ of $X$ and the function $\varphi: (X-K)\to I$ is a proper pluriharmonic function onto  an open interval $I$ of $\mathbf{R}$. Let $\pi: \hat{X}\to X$ be a holomorphic covering such that $\pi^*\mathbf{\xi}=dF$ for some holomorphic function $F$ on $\hat{X}$. Let $h=Re\,F$. Then we have one of the following
two cases:\\
	(a) there exists a holomorphic mapping $p$ from $X$ to a Riemann surface $S$ with connected fibers and a holomorphic one-form $\tau$ on $S$ such that $\mathbf{\xi} =p^*\tau$.\\
	(b) $h^{-1}(r)$ (resp. $F^{-1}(z)$) is connected for any $r\in \mathbf{R}$ (resp.$z\in\mathbf{C}$)  and $\pi_1(h^{-1}(r))\to\pi_1(\hat{X})$ (resp. $\pi_1(F^{-1}(z))\to\pi_1(\hat{X})$) is surjective.\\
\end{theorem}

Note that Simpson's theorem can be viewed as the case when $K=X$ in the above $Theorem $ $\ref{thm1}$ (now we need to drop the assumption that $X$ is non-compact).

Note also that if a complex manifold satisfies the condition in the above theorem, then it automatically becomes a weakly 1-complete manifold. One can see this as follows. Let $I=(a,b)$. We can find a compact subset $K'$ of $X$ such that $K\subset K'$ and the minimal value of $\varphi$ on $\partial K'$ is a real number $c$. Then $\psi:=max(c, \varphi)$ is a continuous plurisubharmonic function on $X$.
 If $b=+\infty$, then $\psi$ is a continuous plurisubharmonic exhaustive function on $X$.
If $b<+\infty$, then $-\log |b-\psi|$ is a continuous plurisubharmonic exhaustive function on $X$.

The proof of  $Theorem $ $\ref{thm1}$ is organized as follows. First we introduce the non-compact Stein Factorization in Section 2. Then in Section 3, for each path-connected component of $A$ (the zero set of $\mathbf{\xi}$), we construct a tubular neighborhood of it (and we call the path-connected component of  $A$ the central fiber of its tubular neighborhood). Now for each tubular neighborhood we construct above, we find a holomorphic function defined on it such that the value of this function equals zero on the central fiber. If every central fiber intersects the boundary of its tubular neighborhood, then we show that the nearby fiber is also connected by using the non-compact Stein factorization. In this case, we will see that $h^{-1}(r)$ (resp. $F^{-1}(z)$) is connected for any $r\in \mathbf{R}$ (resp. $z\in\mathbf{C}$).

  While if there exists a central fiber which does not intersect the boundary of the tubular neighborhood, we will show that there exists a holomorphic mapping from $X$ to a Riemann surface $S$ such that the meromorphic function is constant on every fiber of this holomorphic function mapping. Moreover, the holomorphic one-form $\mathbf{\xi}$ is the pull-back of a holomorphic one-form $\tau$ on $S$.

Before proceeding to the proof, let us first make some simplification for our current situation. 
Recall that (see \cite{r16}) an end of $X$ is an element of 
$$\underleftarrow{\lim}\;\pi_0(X-K)$$
where the limit is taken as $K$ ranges over the compact subsets of $X$.

Let 
$$\eta: \widetilde{X}\to X$$ 
be a universal covering of $X$.
 A filtered end of $X$ is an element of 
$$\underleftarrow{\lim}\; \pi_0[\eta^{-1}(X-K)]$$
where the limit is taken as $K$ ranges over the compact subsets of $X$.
The number of ends of $X$ will be denoted by $e(X)$ and the number of filtered ends of $X$ will be denoted by $\tilde{e}(X)$. We have $\tilde{e}(X)\geqslant k$ if and only if there exists an end decomposition $X-K= E_1\cup...\cup E_m$ such that 
$$\sum_{j=1}^m[\pi_1(X):\Gamma_j]\geqslant k,$$
where $\Gamma_j:= im\,(\pi_1(E_j)\to\pi_1(X))$.

 For the case when $X$ is an open K\"ahler manifold with the number of filtered ends greater than or equal to three, Napier and Ramachandran showed in \cite{r16} that if $X$ has bounded geometry, or is hyperbolic, or is weakly 1-complete, then it admits a proper holomorphic mapping to a Riemann surface. So in this paper we just need to consider the following three cases:
(1) $\tilde{e}(X)=e(X)=1$, \quad(2) $\tilde{e}(X)=e(X)=2$, \quad(3) $\tilde{e}(X)=2, e(X)=1$.

	\section{Non-compact Stein Factorization}
	$\mathbf{Definition}$ Let $f: M\to N$ be a surjective smooth mapping between two smooth manifolds. Let $F$ be a smooth manifold. If for every $x\in N$ there is a neighborhood $U$ of $x$ in $N$ such that there is a diffeomorphism
	$$ \phi: f^{-1}(U)\to U\times F $$
	and the diagram
	$$\xymatrix{
		f^{-1}(U) \ar[rr]^{\phi}\ar[dr]_{f} &   & U\times F \ar[dl]^{pr_1} \\
		& U &}$$
	commutes (here $pr_1$ is the projection to the first factor), then we call the triple $(M,N,f)$ a smooth fiber bundle and $F$ the fiber of $f$.
	
		If the $M$ and $N$ are complex manifolds and the $f$ is holomorphic, we have the following theorem:
	\begin{theorem}[Non-compact Stein Factorization]\label{t3}
		Let $(M, N, f)$ be a smooth fiber bundle, where $M$ and $N$ are two complex manifolds and $f: M\to N$ is a holomorphic submersion. If every fiber of $f$ has only finitely many connected components, then one can factor $f$ into $\phi\circ g$, where the holomorphic mapping $g: M\to N'$ has connected fibers and $\phi: N'\to N$ is a finite holomorphic covering. 
	\end{theorem}
	\textbf{proof}. First we introduce an equivalence relation "$\sim$" on $M$. For any two points $x$ and $y$ in $M$,  we define $x\sim y$ if and only if $x$ and $y$ are contained in the same connected component of a fiber of $f$. We then give $(M/\sim)$ the quotient topology.
	Let 
	$$g:M\to (M/\sim) $$ 
	be the quotient map. For any two distinct points $\bar{x}$ and $\bar{y}$ in $(M/\sim)$, if $$f(g^{-1}(\bar{x}))=z_1\neq z_2=f(g^{-1}(\bar{y})),$$ 
	then let $U_{z_1}$ and $U_{z_2}$ be two disjoint neighborhoods of $z_1$ and $z_2$ in $N$ respectively such that 
	$$f^{-1}(U_{z_1})\cong f^{-1}({z_1})\times U_{z_1}$$
	and 
	$$f^{-1}(U_{z_2})\cong f^{-1}({z_2})\times U_{z_2}.$$ 
	We have that $g(f^{-1}(U_{z_1}))$ and $g(f^{-1}(U_{z_1}))$ are two disjoint neighborhoods of $\bar{x}$ and $\bar{y}$ in $(M/\sim)$ respectively. If $$f(g^{-1}(\bar{x}))=f(g^{-1}(\bar{y}))=z,$$
	then take a  small neighborhood $U_z$ of $z$ in $N$ such that 
	$$ f^{-1}(U_z)\cong f^{-1}(z)\times U_z.$$
	Let 
	$$f^{-1}(z)=C_1\cup C_2\cup...\cup C_m$$ where $C_1, C_2,..., C_m$ are the connected components of $f^{-1}(z)$. Then $g^{-1}(\bar{x})= C_i$ and $g^{-1}(\bar{y})= C_j$ for some $i\neq j$. Thus $g(C_i\times U_z)$ and $g(C_j\times U_z)$ are two disjoint neighborhoods of $\bar{x}$ and $\bar{y}$ in $(M/\sim)$ respectively. In both case $z_1$ and $z_2$ can be separated by two disjoint neighborhoods, so  $(M/\sim)$ is a Hausdorff space. Since 
	$$g(C_i\times U_z)\cong U_z$$
	we know $(M/\sim)$ admits a complex manifold structure. Define 
	$$\phi : (M/\sim)\to N$$ 
	by $\phi(\bar{x})=z$. We have 
	$$\phi^{-1}(U_z)=\displaystyle \mathop{\cup}_{j=1}^{m} g(C_j\times U_z).$$ 
	This implies $\phi$ is an $m$-sheet covering map. In fact it is a holomorphic covering map. If we denote $N'$ the quotient space $(M/\sim)$, we have the following commutative diagram:
	$$\xymatrix{
		M \ar[rr]^{g}\ar[dr]_{f} &   & N' \ar[dl]^{\phi} \\
		& N &}$$
	Locally $\phi^{-1}$ is holomorphic, so $g=\phi^{-1}\circ f$ is also holomorphic.
	$\square$
	\section{Proof of  Theorem 1}

	Let 
	$$A=\{x\in X|\mathbf{\xi}_x=0\}.$$
	 It is an analytic subspace of $X$. Note that the leaves of the holomorphic foliation determined by $\mathbf{\xi}$ lie in the fibers of $\varphi$ outside $K$ and only finitely many components of $A$ can meet the boundary of a relatively compact neighborhood of $K$, so all the components of $A$ are compact.
	 
	 Let $g$ be  a function on $X$ defined by
	$$\begin{matrix}
	g:X\to \mathbf{R}\\
	    \quad\qquad x\mapsto \lvert\mathbf{\xi}_x\vert_.^2 \\
		\end{matrix}$$
    Here $\lvert\mathbf{\xi}_x\vert^2$ is the square of the norm of $\mathbf{\xi}_x$ which is induced by the Hermitian metric on the holomorphic cotangent bundle of $X$. So $g$ is a smooth function and thus the set of its critical values has measure zero by Sard's theorem. This implies that for any $\epsilon_1>0$ there exists an $\epsilon$ in $(0,\epsilon_1)$  such that $g^{-1}(\epsilon)$ is a smooth manifold. Hence for any such $\epsilon$, $g^{-1}(-\infty, \epsilon]$ is a smooth submanifold of $X$ with a smooth boundary $g^{-1}(\epsilon)$.	For such an $\epsilon$, let $$U=g^{-1}(-\infty, \epsilon].$$ 
    Let $A_k $ $(k=1, 2, 3...)$ be the path-connected components of $A$. 
    Define $U_k$ as the connected component of $U$ that contains $A_k$. 
    
    We claim that for any finite many path-connected components of $A$ say $A_1, A_2, ...,A_n$, we can make the $\epsilon$ (depending on $n$) so small that all the $U_k$ $($k=1,...,n$)$ are compact and $U_k\neq U_j$ for $k\neq j, 1\leqslant k, j\leqslant n$. We can also assume that $U_k$ contains no other components of $A$ except $A_k$ for each $k$ ($k=1,...,n$), since $\partial\varphi=\mathbf{\xi}$ on $(X-K)$ and $\varphi$ is a proper function on $(X-K)$.
    
    Suppose we have fixed an $n$, then for each $A_k$ there exist an open neighborhood $V'_k$ of $A_k$ and a holomorphic function $f_k$ defined on $V'_k$ such that 
    $$f_k|_{A_k}=0$$ 
    and 
    $$df_k=\mathbf{\xi}|_{V'_k}.$$
    We can take the above $\epsilon$ such that $\partial U_k\cap f_k^{-1}(0)$ (might be empty) is a smooth submanifold of $\partial U_k$ for every $k$. Here $\partial U_k$ denotes the boundary of $U_k$. Let 
    $$\Delta(0,\delta)=\{z\in\mathbf{C}| |z|<\delta\}.$$ 
    And let $\overline{\Delta(0,\delta)}$  be its closure in the complex plane. Since $\mathop{\cup}_{k=1}^{n}A_k$ has only finitely many path-connected components, we can choose a fixed real number $\delta>0$ for all $k$ such that\\
     i) If $\partial U_k\cap f_k^{-1}(0)$ is not empty, then  the $\delta$ is so small that $\partial U_k\cap f_k^{-1}(v)$ is also smooth for all $v\in\Delta(0,\delta)$ (this is equivalent to that $\partial U_k$ intersects $f_k^{-1}(v)$ transversely for all $v\in\Delta(0,\delta)$).\\
     ii) If $\partial U_k\cap f_k^{-1}(0)$ is  empty, we choose the $\delta$ so small that $f_k^{-1}(\overline{\Delta(0,\delta)})$ is contained in $U_k$. Since $A_k$ is compact, we may also assume that the $\delta$ is small enough that $f_k^{-1}(\overline{\Delta(0,\delta)})$, which is a closed subset of $V'_k$ , is also a closed subset of $X$ and that $f_k$ is surjective when it is restricted to $f_k^{-1}(\overline{\Delta(0,\delta)})$. 
     
     We then define $V_k$ to be $f_k^{-1}(\overline{\Delta(0,\delta)})$. 
    Let $\partial V_k$ denote the boundary of $V_k$. Note that 
    $$\partial V_k=f_k^{-1}(\{z\in\mathbf{C}| |z|=\delta\}).$$ 
    Also let $Int(U_k)$ and $Int(V_k)$ denote the interior of $U_k$ and $V_k$, respectively. Put 
    $$W_k=U_k\cap V_k. $$ 
    We have the following:
    	\begin{lemma}[\cite{r11}]\label{l1} 
    	$W_k$ is a smooth manifold with corners. Its boundary is
    	$$\partial W_k=T_k\cup R_k\cup M_k$$
    	where $T_k=Int(V_k)\cap\partial U_k$, $R_k=\partial V_k\cap Int(U_k)$ and $M_k=\partial V_k\cap \partial U_k$. \\
    	If $\partial U_k\cap f_k^{-1}(0)$ is not empty, then $T_k$ and $R_k$ are smooth pieces and $M_k$ is a smooth corner.\\
    	If $\partial U_k\cap f_k^{-1}(0)$ is empty we have  $T_k$ and $M_k$ are both empty and $\partial W_k=R_k=\partial V_k$ by definition.
    	
    \end{lemma}  
   \textbf{proof}. We just need to prove the case when $\partial U_k\cap f_k^{-1}(0)$ is not empty. Since $ U_k$, as a path-connected component of a smooth manifold with boundary, is smooth we know $\partial U_k$ is smooth. Thus $T_k$ is smooth.
    $R_k$ is smooth because every point in $\partial V_k$ is not a critical point of $f_k$.
    
    We have chosen the 
    $\delta$ so that $\partial U_k\cap f_k^{-1}(v)$ is smooth for all $v\in\Delta(0,\delta)$. This smoothness is equivalent to that $\partial U_k$ and $\partial V_k$ intersect each other transversely, which implies that $M_k$ is smooth. $\square$   
    
    The mapping $f_k|_{W_k}$ is a proper mapping because $W_k$ is a compact set (it is a closed subset of the compact set $U_k$).
     Let
     $$\mathring{W_k}=(W_k-f_k^{-1}(0)-(R_k\cup M_k))=(W_k-f_k^{-1}(0)-\partial V_k) $$ 
     and 
     $$ \mathring{\Delta}(0,\delta)=(\Delta(0,\delta)-\{0\}).$$
      We have
      \begin{lemma}
       The restriction of $f_k$ to $\mathring{W_k}$, 
    $$f_k: \mathring{W_k}\to \mathring{\Delta}(0,\delta)$$ 
    is a fibration. 
    \end{lemma}
    \textbf{proof}.  For any compact $K\subseteq \mathring{\Delta}(0,\delta)$ we have $$(f_k|_{\mathring{W_k}})^{-1}(K)=(f_k|_{W_k})^{-1}(K)$$ is a compact set. So $f_k|_{\mathring{W_k}}$ is also a proper mapping. Let $Int(\mathring{W_k})$ be the interior of $\mathring{W_k}$ and $\partial{\mathring{W_k}}$ the boundary of $\mathring{W_k}$. We know
    $$Int(\mathring{W_k})=\mathring{W_k}-T_k=\mathring{W_k}-\partial U_k$$
    and
    $$\partial{\mathring{W_k}}=\mathring{W_k}\cap T_k=\mathring{W_k}\cap\partial U_k$$
    by $Lemma $ $ \ref{l1}$ above.
      Since every point in $Int(\mathring{W_k})$ is not a critical point of $f_k$ we have $f_k$ is a submersion on $Int(\mathring{W_k})$ . When $\partial U_k\cap f_k^{-1}(0)$ is not empty we have $\partial U_k\cap f_k^{-1}(z)$ is smooth for all $z\in\Delta(0,\delta)$. This implies that  $f_k$ is also a submersion when restricted to $\partial{\mathring{W_k}}=\mathring{W_k}\cap\partial U_k$.  If $\delta$ is sufficiently small $f_k|_{\mathring{W_k}}$ will also be surjective. Therefore, $f_k: \mathring{W_k}\to \mathring{\Delta}(0,\delta)$ is a fiber bundle by Ehresmann's fibration theorem. $\square$
      
    Let $$\widetilde{W_k}=\mathring{W_k}-\partial{\mathring{W_k}}$$ and $\tilde{f_k}$ the restriction of $f_k$ to $\widetilde{W_k}$. Then 
    $$\tilde{f_k}: \widetilde{W_k}\to \mathring{\Delta}(0,\delta)$$ 
    is still a fiber bundle. Since a fiber of $f_k$ is a compact submanifold of $\mathring{W_k}$, it has finitely many connected components. Therefore a fiber of $\tilde{f_k}$ also has finitely many connected components. In fact, if $\tilde{C}$ is a connected component of $f_k^{-1}(z)$, then $C=\tilde{C}-\partial\tilde{C}$ is a connected component of  $\tilde{f_k}^{-1}(v)$. So $\tilde{f_k}^{-1}(z)$ has the same number of connected components as $f_k^{-1}(z)$. Note that here $\partial\tilde{C}$ is the boundary of $\tilde{C}$ and it might be empty.\\
       From  $Theorem $ $\ref{t3}$ and the above lemmas we can easily get the following:
       \begin{corollary}\label{l3} 
       	Consider the fiber bundle
       	 $$\tilde{f_k}: \widetilde{W_k}\to \mathring{\Delta}(0,\delta).$$   
       	We have the following commutative diagram:
       	$$\xymatrix{
       			\widetilde{W_k} \ar[rr]^{g_k}\ar[dr]_{\tilde{f_k}} &   &  \mathring{\Delta}(0,\delta') \ar[dl]^{\phi_k} \\
       		& \mathring{\Delta}(0,\delta) &}$$
       	where $\phi_k(z)=z^m$ and $m$ is the number of connected components of every fiber of $\tilde{f_k}$. Therefore
       	$$\tilde{f_k}=\phi_k\circ g_k=g_k^m.$$
       	\end{corollary}
       \begin{corollary}
       If $\partial U_k\cap f_k^{-1}(0)$ is not empty, then every fiber of $\tilde{f_k}$ has just one connected component.
       \end{corollary}
        \textbf{proof}. by $Corollary $ $ \ref{l3}$ above we have 
        $$d\tilde{f_k}=mg_k^{m-1}dg_k.$$ 
        Both $f_k$ and $g_k$ have continuous extensions to the whole $W_k$, and by Riemann removable singularity theorem the extensions are holomorphic. We then get
        $$df_k=mg_k^{m-1}dg_k.$$
       If $m>1$, then for a point $p$ in $f_k^{-1}(0)=g_k^{-1}(0)$	we have $$df_k|_p=mg_k^{m-1}(p)dg_k|_p=0,$$ 
       and thus $p\in A_k$. This contradicts that
        $\partial U_k\cap f_k^{-1}(0)$ is a smooth submanifold of $\partial U_k$. So we must have $m=1$. $\square$
        
\subsection{Case (b) of theorem 1}

        Let 
    $$\pi: \hat{X}\to X$$
     be a covering space such that 
     $$\pi^*\mathbf{\xi}=dF$$
      for some holomorphic function $F$ on $\hat{X}$. Define $h=Re F$. 
      
    Let 
    $$\psi: \widetilde{X}\to \hat{X}$$ 
    be the universal covering manifold of $\hat{X}$ and $\tilde{h}=\psi^*h$. Of course 
    $$\eta=\pi\circ\psi:\widetilde{X}\to X$$
     is also the universal covering of $X$. We know 
     the set of critical points of $\tilde{h}$ is $\eta^{-1}(A)$. 
     
     If $A$ is empty, then we will have
     $$\hat{X}=h^{-1}(r)\times \mathbf{R}$$
     as a $C^{\infty}$ manifold. In this case, it is trivial  that $h^{-1}(r)$ is connected for all $r\in \mathbf{R}$  and $\pi_1(h^{-1}(r))\to\pi_1(\hat{X})$ is surjective.
     
     Suppose now $A$ is not empty. 
     
   \textbf {subcase (1)  $\tilde{e}(X)=e(X)=1$}
     
     We know in this case for any end decomposition $X-K= E$, we have the homomorphism $\pi_1(E)\to\pi_1(X)$ is onto. This implies that $\eta^{-1}(E)$ is connected. So we can assume $\eta^*\varphi=2 \tilde{h}$ on $\eta^{-1}(E)$. We may also assume the boundary of $K$ is given by $\varphi^{-1}(s)$ for some regular value of $\varphi$.

      Suppose that for every $k$, $\partial U_k\cap f_k^{-1}(0)$ is not empty. (If  $\partial U_k\cap f_k^{-1}(0)$ is empty for some $k$, then we will show there exists a holomorphic mapping from $X$ to a Riemann surface in subsection 3.2.) We will  prove that,
       \begin{lemma}\label{l4}
        If for every $k$ the set $\partial U_k\cap f_k^{-1}(0)$ is not empty, then
        $\tilde{h}^{-1}[r,+\infty)$ is connected for every $r\in\mathbf{R}$.
        \end{lemma}
     \textbf{proof}. 
     Suppose that $\tilde{h}^{-1}[r,+\infty)$ is not connected and $B_0$ and $B_1$ are two different connected components of $\tilde{h}^{-1}[r,+\infty)$.
     Now let 
     $$\gamma: [0,1]\to \widetilde{X}$$
     be a curve in $\widetilde{X}$ such that $\gamma(0)\in B_0$ and $\gamma(1)\in B_1$. If we define $r'$ to be 
     $$r'=\min \limits_{0\leqslant t\leqslant 1} \tilde{h}(\gamma(t)),$$
     we can then view $\gamma$ as a curve in 
     $$\tilde{h}^{-1}[r',+\infty).$$
     Since the image of $\tilde{h}^{-1}[r',r]$ under $\eta$ is contained in $K\cup\varphi^{-1}[s,\frac{1}{2}r]$ and there are only finite many components of $A$ contained in $K\cup\varphi^{-1}[s,\frac{1}{2}r]$ (Note that $\varphi$ is proper so $\varphi^{-1}[s,\frac{1}{2}r]$ is compact) we can assume these components are $A_1, A_2,...,A_n$ and fix this $n$.
     
     We know $$\displaystyle \mathop{\cup}_{k=1}^{n}\eta^{-1}(W_k)$$
     is a neighborhood of $\eta^{-1}(\mathop{\cup}_{k=1}^{n}A_k)$. 
     
     Let $c$ be a real number such that $0<5c<\delta$ and $3c$ is a regular value of $f_k$ for $k=1,...,n$. For each $k \;(1\leqslant k\leqslant n)$, we can find two open neighborhoods $W''_k$ and $W'_k$  of $A\cap W_k$ such that 
     
     $$W''_k\Subset W'_k\Subset W_k\cap f_k^{-1}(\overline{\Delta(0,c)})$$
    and 
    $$\lvert\mathbf{\xi}_x\vert^2 \geqslant \mu$$
     for a real $\mu>0$ when $x$ is not contained in  $W'_k$. Then we can construct a vector field $\textbf{v}$ on $X$ such that 
     $$\textbf{v}_x=0$$
     for $x$ contained in the closure of $W''_k$, and
     $$(Re\mathbf{\xi})(\textbf{v}_x)=1$$ 
     for $x$ not in
      $W'_1\cup W'_2\cup...\cup W'_n,$
      and $\textbf{v}_x$ is tangent to $T_k$ when $x\in T_k$. It is possible to choose such a vector field since $\lvert\mathbf{\xi}_x\vert^2 \geqslant \mu$ when $x$ is not in
      $W'_1\cup W'_2\cup...\cup W'_n$ and $T_k$ intersects $f_k^{-1}(v)$ transversely for all $v\in\Delta(0,\delta)$. Lift this vector field to a vector field $\tilde{\textbf{v}}$ on $\widetilde{X} $. Then note that $d\tilde{h}=Re(\eta^*\mathbf{\xi})$, we have
      $$d\tilde{h}(\tilde{\textbf{v}}_y)=1$$ 
      when $y$ is not in $$\displaystyle \mathop{\cup}_{k=1}^{n}\eta^{-1}(W'_k),$$
    which is a neighborhood of $\eta^{-1}(\mathop{\cup}_{k=1}^{n}A_k)$. And we also have that $\tilde{\textbf{v}}_y$ is tangent to $\eta^{-1}(T_k)$ when $y\in\eta^{-1}(T_k)$. Since $\widetilde{X}$ is complete and $\tilde{\textbf{v}}$ has bounded length, the flow $\phi_t$ generated by $\tilde{\textbf{v}}$ exists for all $t$. 
     
     Let $W^*$ be a path-connected component of $\eta^{-1}(W_k)$. For $z\in  \mathring{\Delta}(0,\delta)$, by $corollary$ $2$ above and a result of Nori (\textbf{Lemma 1.5} of \cite{r9}), we know 
    $$\pi_1(f_k^{-1}(z))\to \pi_1(W_k)$$
    is surjective. So the image of $\pi_1(f_k^{-1}(z))$ meets every coset of the subgroup $\eta_*\pi_1(W^*)$.  Then by a theorem in algebraic topology (for example, \textbf{Proposition 11.2} of Chapter Five in \cite{r6}) we know that $\eta^{-1}(f_k^{-1}(z))\cap W^*$ is connected. 
    
     If $\gamma(0)$ is contained in some $\eta^{-1}(W_k)$ and $f_k(\eta(\gamma(0)))\neq0$, we use a path contained in some  $\eta^{-1}(f_k^{-1}(z))$ to join $\gamma(0)$ and a point in $\eta^{-1}(T_k)$. We then get a new curve which we still denote by $\gamma$ and we let $\gamma(0)$ be the point in $\eta^{-1}(T_k)$. Because $d\tilde{h}-Re\,df_k=0$, the new curve is still in $\tilde{h}^{-1}[r',+\infty)$. If  $f_k(\eta(\gamma(0)))=0$, then since $f_k^{-1}(0)$ is connected and $f_k^{-1}(0)\cap T_k$ is not empty we can also get a curve with $\gamma(0)\in\eta^{-1}(T_k)$. We do the same thing for $\gamma(1)$. If both  $\gamma(0)$ and  $\gamma(1)$ are not in any $\eta^{-1}(W_k)$ we do nothing.
     
     Since the image of $\gamma$ is compact, it only intersects finitely many path-connected components of each 
     $$\eta^{-1}(f_k^{-1}(\overline{\Delta(0,3c)})\cap W_k).$$ 
     For the parts of $\gamma$ that are outside 
     $$\eta^{-1}(f_k^{-1}(\Delta(0,3c))\cap W_k),$$ 
     we may assume all the end points of these parts are contained in $\eta^{-1}(T_k)$. (If not, we can join the end point to some point in $\eta^{-1}(T_k)$ by a path that is contained in some 
     $$\eta^{-1}(f_k^{-1}(z))\cap W^*.)$$
     Note that $|z|=3c$ in this case. We then get a new path. Because 
     $$d\tilde{h}-Re\,df_k=0,$$ 
     the new path is still in $\tilde{h}^{-1}[r',+\infty)$.
     
      Then we can move these paths by $\phi_c$ into $\tilde{h}^{-1}[r'+c,+\infty)$ since $W'_k\Subset  f_k^{-1}(\overline{\Delta(0,c)})$ and $d\tilde{h}(\tilde{\textbf{v}}_y)=1$ when $y$ is not in $\displaystyle \mathop{\cup}_{k=1}^{n}\eta^{-1}(W'_k).$ (Remember that $\phi_t$ is the flow generated by $\tilde{\textbf{v}}$.) The end points of these paths will be still contained in $\eta^{-1}(T_k)$, because 
     $\tilde{\textbf{v}}$ is tangent to $\eta^{-1}(T_k)$ and $5c<\delta$. So in fact the end points of these paths are all contained in
      $$\tilde{h}^{-1}[r'+c,+\infty)\cap\eta^{-1}(T_k).$$ 
     
      For two points $y'$ and $y''$  in the same piece of some $\eta^{-1}(T_k)$, we can find a path that is contained in that piece to join them. If $y'$ and $y''$ are contained in $\tilde{h}^{-1}[r'+c,+\infty),$ this path can also be chosen to be contained in $\tilde{h}^{-1}[r'+c,+\infty).$ For a path-connected component $W^*$ of $\eta^{-1}(W_k)$ we can join all the pieces of $W^*\cap \eta^{-1} (T_k)$ by some $\eta^{-1}(f_k^{-1}(z))\cap W^*.$ If $W^*\cap \eta^{-1} (T_k)$ contains some point in 
       $\tilde{h}^{-1}[r'+c,+\infty)$ we can choose the $\eta^{-1}(f_k^{-1}(z))\cap W^*$ to be contained in $\tilde{h}^{-1}[r'+c,+\infty)$ too. This implies all the paths contained in $\tilde{h}^{-1}[r'+c,+\infty)$ which we construct above can be jointed together by paths contained in $\tilde{h}^{-1}[r'+c,+\infty)$. Thus $B_0$ and $B_1$ can be jointed by a curve in $\tilde{h}^{-1}[r'+c,+\infty)$. Since $c$ is a fixed positive number, by repeating what we did above we can finally get a curve contained in $\tilde{h}^{-1}[r,+\infty)$ that joins $B_0$ and $B_1$. This is a contradiction. So $\tilde{h}^{-1}[r,+\infty)$ must be connected. $\square$
       
      Similarly, we can prove that $\tilde{h}^{-1}(-\infty,r]$ is also connected. 
      
      Since the function $\tilde{h}$ is real analytic, we know both $\tilde{h}^{-1}[r,+\infty)$ and $\tilde{h}^{-1}(-\infty,r]$ are triangulable and $\tilde{h}^{-1}(r)$ can be realized as a subcomplex of each of the triangulations (see \cite{r5}). Consider the Mayer-Vietoris sequence of reduced homology with real coefficient
    $$H_1(\widetilde{X},\mathbf{R})\to H_0^\#(\tilde{h}^{-1}(r),\mathbf{R})\to H_0^\#(\tilde{h}^{-1}(-\infty,r],\mathbf{R})\oplus H_0^\#(\tilde{h}^{-1}[r,+\infty),\mathbf{R}).$$
    We know 
    $$H_0^\#(\tilde{h}^{-1}(-\infty,r],\mathbf{R})= H_0^\#(\tilde{h}^{-1}[r,+\infty),\mathbf{R})=0$$
     since $\tilde{h}^{-1}[r,+\infty)$ and $\tilde{h}^{-1}(-\infty,r]$ are connected. We also have 
     $$H_1(\widetilde{X},\mathbf{R})=0$$
      since $\widetilde{X}$ is the universal covering space of $\hat{X}$. This implies  $$H_0^\#(\tilde{h}^{-1}(r),\mathbf{R})=0$$
       and thus $\tilde{h}^{-1}(r)$ is connected. Therefore $h^{-1}(r)$ is also connected because $h^{-1}(r)=\eta(\tilde{h}^{-1}(r))$.
       
       Let 
       $$i:h^{-1}(r)\hookrightarrow\hat{X}$$
       be the inclusion. Let 
       $$\Gamma=i_*\pi_1(h^{-1}(r))$$
       be the image of $\pi_1(h^{-1}(r))$ in $\pi_1(\hat{X}).$ We have the following diagram:
       $$\begin{tikzcd}
        &&&& \widetilde{X}\arrow{ddd}[swap]{\psi_1} \arrow[bend left=30]{dddddd}{\psi}\\
        &&&&\\
        &&&&\\
        &&&&_{\Gamma}\backslash^{\widetilde{X}} \arrow{ddd}[swap]{\psi_2}\\
        &&&&\\
        &&&&\\
       h^{-1}(r)\arrow{rrrr}{i}\arrow{uuurrrr}{i'}&&&&\hat{X} 
       \end{tikzcd}  $$
       where $i'$ is the lifting of $i$ to ${\Gamma}\backslash{\widetilde{X}}$ and $\psi$ factors into the composite of $\psi_1$ and $\psi_2$.
       
        If $\Gamma$ is not equal to $\pi_1(\hat{X})$, then it can not meet every coset of ${\psi_2}_*\pi_1({\Gamma}\backslash{\widetilde{X}})$ in $\pi_1(\hat{X})$. So in this case $i'(h^{-1}(r))$ is not connected (see \textbf{Proposition 11.2} of Chapter Five in \cite{r6}). But $i'(h^{-1}(r))$ equals $\psi_1(\tilde{h}^{-1}(r))$ which is a continuous image of a connected set. This is a contradiction. Thus we must have $\Gamma=\pi_1(\hat{X})$. That is,
       $$ \pi_1(h^{-1}(r))\to\pi_1(\hat{X})$$
       is surjective.
       
Now consider the complex case. Construct another vector field $\textbf{u}$ on $X$ such that $$\textbf{u}_x=0$$
when $x$ is contained in the closure of $W''_k$, and
$$(Im\,\mathbf{\xi})(\textbf{u}_x)=1$$
$$(Re\,\mathbf{\xi})(\textbf{u}_x)=0,$$
when $x$ is not in
$W'_1\cup W'_2\cup...\cup W'_n,$
and $\textbf{u}_x$ is tangent to $T_k$ when $x\in T_k$.   Lift this vector field to a vector field $\tilde{\textbf{u}}$ on $\widetilde{X}.$  Let the flow generated by $\tilde{\textbf{u}}$ be $\psi_t$. For every $r\in \mathbf{R}$, we know $\tilde{h}^{-1}(r)$ is a real analytic subspace of $\widetilde{X}$. 
Let us denote the set of singular points of $\tilde{h}^{-1}(r)$ by $$Sing(\tilde{h}^{-1}(r)).$$  
We have 
$$Sing(\tilde{h}^{-1}(r))\subseteq\tilde{h}^{-1}(r)\cap\eta^{-1}(A).$$
Let 
$$Reg(\tilde{h}^{-1}(r))=\tilde{h}^{-1}(r)-Sing(\tilde{h}^{-1}(r)).$$
By the construction of $\tilde{\textbf{u}}$, we can view the restriction of $\tilde{\textbf{u}}$ on 
$Reg(\tilde{h}^{-1}(r))$ 
as a vector field on 
$Reg(\tilde{h}^{-1}(r))$. Define $l=ImF$ and $\tilde{l}=\psi^*l$. Denote the restriction of $\tilde{l}$ to $\tilde{h}^{-1}(r)$  by $\tilde{l}_r$. We have 
$$d\tilde{l}_r(\tilde{\textbf{u}}) =1$$ 
on $Reg(\tilde{h}^{-1}(r))$ for every $r$.

We can also see that the flow generated by the restriction of $\tilde{\textbf{u}}$ on $Reg(\tilde{h}^{-1}(r))$ is just the restriction of $\psi_t$ on $Reg(\tilde{h}^{-1}(r))$, so we use the same symbols $\tilde{\textbf{u}}$ and $\psi_t$ to denote the restrictions of $\tilde{\textbf{u}}$ and $\psi_t$ on $Reg(\tilde{h}^{-1}(r))$. The $\tilde{\textbf{u}}$ on $Reg(\tilde{h}^{-1}(r))$ is a complete vector field since $\psi_t$ exists for all $t$.  

Let 
$$\zeta :H_r\to\tilde{h}^{-1}(r)$$ 
be the universal covering space of $\tilde{h}^{-1}(r)$. $H_r$ has a real analytic space structure and the set of singular points of $H_r$ is $\zeta^{-1}(Sing(\tilde{h}^{-1}(r))).$ That is,
$$Sing(H_r)=\zeta^{-1}(Sing(\tilde{h}^{-1}(r))).$$
Also, we denote 
$$Reg(H_r)=H_r-Sing(H_r).$$
Lift $\tilde{\textbf{u}}$ to a vector field $\widetilde{\textbf{U}}$ on $Reg(H_r)$. The flow generated by $\widetilde{\textbf{U}}$ is 
$$\Psi_t=\psi_t\circ\zeta.$$ 
It exists for all $t$. Let 
$$\widetilde{L}_r=\tilde{l}_r\circ\zeta.$$
We have the following lemma,
\begin{lemma}
	If for every $k$ the set $\partial U_k\cap f_k^{-1}(0)$ is not empty, then for a fixed $r\in\mathbf{R}$,
	$\widetilde{L}_r^{-1}[s,+\infty)$ is connected for every $s\in\mathbf{R}$.
\end{lemma}
\textbf{proof}. The proof is almost the same as the proof of $Lemma$ $\ref{l4}$  above. Suppose that $\widetilde{L}_r^{-1}[s,+\infty)$ is not connected with $B_0$ and $B_1$ two different path-connected components of $\widetilde{L}_r^{-1}[s,+\infty)$. let 
$$\gamma: [0,1]\to H_r$$
be a curve in $H_r$ such that $\gamma(0)\in B_0$ and $\gamma(1)\in B_1$. If we define $s'$ to be 
$$s'=\min \limits_{0\leqslant t\leqslant 1} \widetilde{L}_r(\gamma(t)),$$
we can then view $\gamma$ as a curve in 
$$\widetilde{L}_r^{-1}[s',+\infty).$$
First we show that for any path-connected component $W^{**}$ of $\zeta^{-1}(\eta^{-1}(W_k))$, the subset
$$\zeta^{-1}(\eta^{-1}(f_k^{-1}(z)))\cap W^{**}$$ is path-connected.
For every $W_k$ we know there exists a constant real number $c_k$ such that 
$$[(Ref_k)\circ\eta]^{-1}(r+c_k)=\tilde{h}^{-1}(r).$$
So in order to show $\zeta^{-1}(\eta^{-1}(f_k^{-1}(z)))\cap W^{**}$ is connected, we just need to show there is a surjective homomorphism:
$$\pi_1(f_k^{-1}(z))\to\pi_1[(Ref_k)^{-1}(r+c_k)\cap W_k].$$
We can show there is a retraction (\textbf{Lemma 10} in \cite{r11})
$$R:W_k\to (Ref_k)^{-1}(r+c_k)\cap W_k .$$ 
This implies the homomorphism 
$$R_*: \pi_1(W_k)\to\pi_1[(Ref_k)^{-1}(r+c_k)\cap W_k]$$
is surjective. As in the proof in $Lemma$ $\ref{l4}$ we know the homomorphism
$$\pi_1((f_k^{-1}(z))\to\pi_1(W_k)$$
is surjective for $z\in  \mathring{\Delta}(0,\delta)$ with $Rez=r+c_k$, so the composition 
$$\pi_1((f_k^{-1}(z))\to\pi_1(W_k)\to\pi_1[(Ref_k)^{-1}(r+c_k)\cap W_k]$$
is surjective. 

Now for the parts of $\gamma$ that lie outside 
$$\displaystyle \mathop{\cup}_{k=1}^{n}\zeta^{-1}\eta^{-1}(W_k)$$
we use the flow $\Psi_t$ to move them into $\widetilde{L}_r^{-1}[s'+c,+\infty)$ (note that in the proof of $Lemma$ $\ref{l4}$, we do not need the flow existing near $\eta^{-1}(\mathop{\cup}_{k=1}^{n}A_k)$, so we can apply the same method here). Then as in the proof of $Lemma$ $\ref{l4}$, we may assume all the intersections of $\gamma$ with $\displaystyle \mathop{\cup}_{k=1}^{n}\zeta^{-1}\eta^{-1}(W_k)$ are contained in 
$$\displaystyle \mathop{\cup}_{k=1}^{n}\zeta^{-1}\eta^{-1}(T_k),$$
so we can use some $\zeta^{-1}(\eta^{-1}(f_k^{-1}(z)))\cap W^{**}$ that contained in $\widetilde{L}_r^{-1}[s'+c,+\infty)$ to join these intersection points. This implies we can find a curve that joins $B_0$ and $B_1$ in $\widetilde{L}_r^{-1}[s'+c,+\infty)$. Since $c$ is a fixed positive number, by repeating this process we can get a curve in $\widetilde{L}_r^{-1}[s,+\infty)$ that joins $B_0$ and $B_1$. But this contradicts that $B_0$ and $B_1$ are two different path-connected components of $\widetilde{L}_r^{-1}[s,+\infty)$.$\square$

Similarly, we can prove $\widetilde{L}_r^{-1}(-\infty,s]$ is connected. Then as in the proof of the real case above, by the Mayer-Vietoris sequence and that $H_r$ is simply connected we get $\widetilde{L}_r^{-1}(s)$ is connected. So $\tilde{l}_r^{-1}(s)$ as the continuous image of $\widetilde{L}_r^{-1}(s)$ is also connected. In a similar way to the real case, we can prove 
$$\pi_1(\tilde{l}_r^{-1}(s))\to\pi_1(h^{-1}(r))$$
is surjective. Then combining the result that $ \pi_1(h^{-1}(r))\to\pi_1(\hat{X})$ is surjective, we get
$$\pi_1(\tilde{l}_r^{-1}(s))\to\pi_1(\hat{X})$$
is surjective. Notice that $r$ and $s$ are arbitrary, finally we have for any $z\in \mathbf{C}$, $F^{-1}(z)$ is connected and 
$$\pi_1(F^{-1}(z))\to \pi_1(\hat{X})$$
is surjective.

\textbf{Remark 1} Note in the reasoning above we just use the fact that only finitely many components of $A$ are contained in $\eta(\tilde{h}^{-1}[r',r])$ for any two arbitrary real numbers $r$ and $r'$.

\textbf{Remark 2} In fact in subcase (1), we did not use the K\"ahler condition. So we may just assume $X$ is a non-compact connected complex manifold in this case.\\

 \textbf {subcase (2) $\tilde{e}(X)=e(X)=2$}

Let $X-K=E_1\cup E_2$ be an end decomposition. 

Let $\alpha$ and $\beta$ be two regular values of $\varphi$ such that $\varphi^{-1}(\alpha)\subset E_1$ and $\varphi^{-1}(\beta)\subset E_2$. We know the closed submanifold of $X$ bounded by $\varphi^{-1}(\alpha)$ and $\varphi^{-1}(\beta)$ is a compact submanifold of $X$ with smooth boundary.  Let $X'$ denote its interior. Since $\varphi$ is a pluriharmonic function, we may assume that there exists a pluriharmonic function $\rho_0$ defined on the neighborhoods $U_{\alpha}$ and $U_{\beta}$ of $\varphi^{-1}(\alpha)$ and $\varphi^{-1}(\beta)$ respectively such that 
$$\varphi^{-1}(\alpha)\cup\varphi^{-1}(\beta)=\rho_0^{-1}(-1)$$ 
and $\rho_0<-1$ in  $(U_{\alpha}\cup U_{\beta})\cap X'$. We extend $\rho_0$ to the whole $X'$ such that $\rho_0<-1$ everywhere on $X'$. Define
$$\rho_1:=max\{\rho_0,-\frac{1}{2}\}.$$
We know $\rho_1$ is a continuous negative plurisubharmonic  function on $X'$. We then define 
$$\rho_2:=-\log(-\rho_1-1).$$ 
We have $\rho_2$ is a continuous exhaustive plurisubharmonic  function on $X'$. By a theorem of (\cite{r19}), given a positive continuous function $\alpha$ on $X'$, there exists a smooth function $\rho_3$ on $X'$ such that
$$\rho_2<\rho_3<\rho_2+\alpha$$ 
and 
$$i\partial\bar{\partial}\rho_3\geqslant-\alpha \hat{\mu},$$
where $\hat{\mu}$ is the K\"ahler form on $X$. Now let 
$$\mu=\hat{\mu}+i\partial\bar{\partial}[(\rho_3)^2]=\hat{\mu}+2i\rho_3\partial\bar{\partial}\rho_3+2i\partial\rho_3\wedge\bar{\partial}\rho_3\geqslant(1-2(\rho_2+\alpha))\hat{\mu}+2i\partial\rho_3\wedge\bar{\partial}\rho_3.$$
By choosing $\alpha>0$ so that $(\rho_2+\alpha)\alpha<\frac{1}{4}$ on $X'$, we get
$$\mu\geqslant\frac{1}{2}\hat{\mu}+2i\partial\rho_3\wedge\bar{\partial}\rho_3.$$
Since $$|d\rho_3|_{\mu}\leqslant\sqrt{2}|\partial\rho_3|_{\mu}\leqslant1,$$ 
we know $\mu$ is complete (because this implies $|\rho_3(x)-\rho_3(y)|\leqslant dist(x, y)$, where $dist(,)$ is the distance function  defined by $\mu$ on $X'$). Now let us fix this K\"ahler metric on $X'$.

 We know $-\rho_1$ is a nonconstant positive continuous superharmonic function, since a plurisuperharmonic function on a K\"ahler manifold is superharmonic. This implies $X'$ is hyperbolic (Recall that a Riemannian manifold is called hyperbolic if and only if  it admits a positive Green's function. See \textbf{Definition 1.1} and \textbf{Characterizations of hyperbolicity} of \cite{r17}).
 
  Let 
  $$E'_1=X'\cap E_1\quad E'_2=X'\cap E_2,$$ 
  and let $G(x,y)$ be the Green's function on $X'$. We call a sequence $\{x_m\}$ approaching $\infty$ in a hyperbolic manifold a regular sequence, if for a fixed $y_0$ the Green's function $G(x_m,y_0)\to0$ as $m\to\infty$. If every sequence approaching $\infty$ is regular, then we say the hyperbolic manifold itself is regular. Since $(2+\rho_1)$ is a plurisubharmonic (and thus subharmonic) function on $X'$ and
  $$\lim\limits_{x \to \infty }(2+\rho_1)|_{E'_1}(x)=1\quad\lim\limits_{x \to \infty }(2+\rho_1)|_{E'_2}(x)=1$$ (remember that $\varphi^{-1}(\alpha)\cup\varphi^{-1}(\beta)=\rho_0^{-1}(-1)=\rho_1^{-1}(-1))$,
  by part (c) of \textbf{Lemma 1.4} of  \cite{r17} and the construction of $E'_1$ and $E'_2$, we know $X'$ is a regular hyperbolic manifold. 
 
 Then by \textbf{Theorem 2.6} of \cite{r17}, there exists a proper pluriharmonic function $\rho$ defined on $X'$ such that $$\lim\limits_{x \to \infty }\rho|_{E'_1}(x)=1 \quad  \lim\limits_{x \to \infty }\rho|_{E'_2}(x)=0.$$ 
 Let $a$ and $b$ be two regular values of $\rho$, we then have 
\begin{align}
& \quad\int_{\rho^{-1}[a, b]}\partial\rho\wedge\overline{\partial\rho}\wedge\mathbf{\xi}\wedge\overline{\mathbf{\xi}}\wedge\mu^{d-2}\\ &=b\int_{\rho^{-1}(b)}\overline{\partial\rho}\wedge\mathbf{\xi}\wedge\overline{\mathbf{\xi}}\wedge\mu^{d-2} -a\int_{\rho^{-1}(a)}\overline{\partial\rho}\wedge\mathbf{\xi}\wedge\overline{\mathbf{\xi}}\wedge\mu^{d-2}\\
&=b\int_{\rho^{-1}(b)}\overline{\partial}(\rho\,\mathbf{\xi}\wedge\overline{\mathbf{\xi}}\wedge\mu^{d-2})-a\int_{\rho^{-1}(a)}\overline{\partial}(\rho\,\mathbf{\xi}\wedge\overline{\mathbf{\xi}}\wedge\mu^{d-2})\\
&=b\int_{\{\text{boundary of }  \rho^{-1}(b)\}}\rho\,\mathbf{\xi}\wedge\overline{\mathbf{\xi}}\wedge\mu^{d-2}-a\int_{\{\text{boundary of }  \rho^{-1}(a)\}}\rho\,\mathbf{\xi}\wedge\overline{\mathbf{\xi}}\wedge\mu^{d-2}\\
&=0
\end{align}
The last equality is because $\rho$ is a smooth proper function, $ \rho^{-1}(a)$ and $ \rho^{-1}(b)$ are thus smooth compact manifolds without boundaries.

Since $a$ and $b$ are arbitrary, this implies  
$$\partial\rho\wedge\mathbf{\xi}=0$$
on $X'$.
Now if $\partial\rho$ and $\mathbf{\xi}$ are linearly independent, then we have a proper holomorphic mapping from $X$ to a Riemann surface (we will deal with this case in section 3.2). Otherwise we have 
$$\partial\rho=c\,\mathbf{\xi}$$
for some nonzero constant $c$.
so we have
 $$\partial\rho=\partial( c\,\varphi)$$
 on $X'-K$.
 We can now assume $c$ is real since both $\rho$ and $\varphi$ are real. Then we have
$$\eta^*(\frac{1}{c}\rho)=2 \tilde{h}.$$
Since $\rho$ is proper and $\alpha$ and $\beta$ are arbitrary, by $ remark \;1$ at the end of subcase (1), we can repeat what we did in subcase (1).\\

  \textbf {subcase (3) $\tilde{e}(X)=2$, $e(X)=1$}
 
 Let $X-K= E$ be an end decomposition. Let $\Gamma:=im\,(\pi_1(E)\to\pi_1(X))$. We have $[\pi_1(X):\Gamma]=2$.  In this case we know there exists a double covering 
 $$\pi': \breve{X}\to X$$ 
 such that $\Gamma=im\,(\pi_1(\breve{X})\to \pi_1(X))$ and $\pi'$ maps some component $\Omega_1$ of $(\pi')^{-1}(E)$  isomorphically onto $E$. In fact, this $\breve{X}$ can be constructed as the quotient ${\Gamma}\backslash{\widetilde{X}}$, where $\widetilde{X}$ is some universal covering space of $X$. 
  Let 
  $$i:E\to X$$ 
  denote be the inclusion map. Since $$im\,(\pi_1(E)\to\pi_1(X))=\Gamma=im\,(\pi_1(\breve{X})\to \pi_1(X))$$
  we know $i$ lifts to a map 
  $$i': E\to \breve{X}.$$
   Let 
   $$\Omega_1=i'(E).$$ 
   Since $[\pi_1(X):\Gamma]=2$, we know there is only one nontrivial covering transformation $$T':\breve{X}\to\breve{X}.$$ 
   Then let 
   $$\Omega_2=T'(i'(E)).$$
   From \textbf{Proposition 11.2} of Chapter Five in \cite{r6} we can see that $(\pi')^{-1}(E)$ is not connected. So $\Omega_1$ and $\Omega_2$ are the two connected components of $(\pi')^{-1}(E)$. Since $(\pi')^{-1}(K)$ is also compact, we have 
  $$\breve{X}-(\pi')^{-1}(K)=\Omega_1\cup \Omega_2$$
  is an end decomposition of $\breve{X}$ and $e(\breve{X})=2$. 
  
  Let $$\eta':\widetilde{X}\to \breve{X}$$
  be a universal covering of $\breve{X}$. Then 
  $$( \eta')^*(\pi')^*\mathbf{\xi}=d\widetilde{G}$$
   for some holomorphic function $ \widetilde{G}$ on $\widetilde{X}$. Let 
   $$\tilde{g}=Re \,\widetilde{G}.$$
   Like in subcase (2), we can prove every $\tilde{g}^{-1}(r)$ is connected and $$\pi_1(\tilde{g}^{-1}(r))\to\pi_1(\widetilde{X})$$
   is onto.  
   
   If $$\pi: \hat{X}\to X$$
   is  a covering space such that 
   $$\pi^*\mathbf{\xi}=dF$$
   for some holomorphic function $F$ on $\hat{X}$, let $h=Re\, F$. \\
   Let 
   $$\psi: \widetilde{X}\to \hat{X}$$ 
   be the universal covering manifold of $\hat{X}$ and $\tilde{h}=\psi^*h$. \\
   We know 
   $$\tilde{g}^{-1}(r)=T (\tilde{h}^{-1}(r))$$
   for some covering transformation
   $$T: \widetilde{X}\to\widetilde{X}.$$
   So we have the same results as in subcase (1).

      \subsection{Case (a) of theorem 1}
    If there exist a $k$ such that $f_k^{-1}(0)\cap U_k=\emptyset$, then there exist a $\delta$ such that for this $k$
    $$f_k^{-1}(\Delta(0,\delta))\cap U_k=\emptyset.$$ 
    So the foliation defined by $\mathbf{\xi}$ has at least one compact leaf.
    By \textbf{Lemma 1.5} of \cite{r15} (and the remarks following it) we know there exists a proper surjective holomorphic mapping $p$ from $X$ to a Riemann surface (See also \cite{r21}). Since the mapping is constant on every leaf of the foliation defined by $\mathbf{\xi}$, we know $\mathbf{\xi} =p^*\tau$ where $\tau$ is a holomorphic one-form.
    
    If in subcase (2) of section 3.1 we have that $\partial\rho$ and $\mathbf{\xi}$ are linearly independent, then together with that $\partial\rho\wedge\mathbf{\xi}=0$, we know $X'$ admits a proper holomorphic mapping onto a Riemann surface by  \textbf{Lemma 1.4} in \cite{r15}. Then by \textbf{Lemma 1.5} of \cite{r15} and its remarks again, we know there exists a proper surjective holomorphic mapping $p$ from $X$ to a Riemann surface. Also we have $\mathbf{\xi} =p^*\tau$ for some holomorphic one-form $\tau$ on that Riemann surface.
    
    If $\tilde{e}(X)\geqslant3$, we know there always exists a proper surjective holomorphic mapping form $X$ to a Riemann surface by \textbf{Theorem 0.1} in \cite{r16}.
    \section{Acknowledgment}
    
    I would like to thank my advisor Professor Ramachandran  who has provided me with a lot of valuable advice on this question. His patient teaching is a great encouragement for me to study math in the future.

\end{document}